\documentclass[12pt,A4]{amsart}

\usepackage{amsmath}
\usepackage{amssymb}

\newtheorem{theorem}{Theorem}[section]

\newtheorem{lemma}[theorem]{Lemma}
\newtheorem{prop}[theorem]{Proposition}
\newtheorem{cor}[theorem]{Corollary}
\theoremstyle{remark}

\newtheorem{remark}[theorem]{Remark}

\def\N{{\mathbb N}}

\def\Q{{\mathbb Q}}
\def\C{{\mathbb C}}
\def\T{{\mathbb T}}

\def\Zp{{{\mathbb Z}_p}}
\def\z{{\mathcal Z}}
\def\Z{{\mathbb Z}}
\def\U{{\mathcal{U}}}
\def\K{{\mathcal{K}}}

\newcommand{\lsp}{\operatorname{span}}
\newcommand{\id}{\operatorname{id}}
\newcommand{\Ind}{\operatorname{Ind}}

\newcommand{\Aut}{\operatorname{Aut}}
\newcommand{\Prim}{\operatorname{Prim}}
\newcommand{\LHS}{\operatorname{LHS}}
\newcommand{\RHS}{\operatorname{RHS}}
\newcommand{\order}{{o}}
\newcommand{\whitesquare}{\hfill $\whitesquare$\newline\vspace{0.4cm}}

\numberwithin{equation}{section}

\begin{document}

\title[Subquotients of Hecke $C^*$-algebras]{\boldmath{Subquotients of Hecke 
$C^*$-algebras}}

\author[N. Brownlowe]{Nathan Brownlowe}
\address{School of Mathematical and Physical Sciences, University of
  Newcastle, NSW 2308, Australia}
\email{Nathan.Brownlowe@studentmail.newcastle.edu.au} 

\author[N.S. Larsen]{Nadia S. Larsen}
\address{Department of Mathematics, Institute for Mathematical Sciences, 
University of Copenhagen, Universitetsparken 5, 2100 Co\-penhagen, 
Denmark}
\address[Current address] 
{Department of Mathematics, University of Oslo, P. O. Box 1053, 
Blindern, N-0316 Oslo, Norway}
\email{nadiasl@math.uio.no}

\author[I.F. Putnam]{Ian F. Putnam}

\address{Department of Mathematics and Statistics, University of 
Victoria, British Colum\-bia V8W 3P4, Canada}
\email{putnam@math.uvic.ca}

\author[I. Raeburn]{Iain Raeburn}
\address{School of Mathematical and Physical Sciences, University of 
  Newcastle, NSW 2308, Australia}
\email{iain@frey.newcastle.edu.au}
\date{September 17, 2003}

\thanks{This research was supported by the Carlsberg Foundation, the Natural 
Sciences and Engineering Research Council of Canada and the Australian 
Research Council.}

\maketitle

In  \cite{BC},  Bost and Connes studied a particular Hecke 
$C^*$-algebra $\mathcal{C}_\Q$ arising in number 
theory. The algebra $\mathcal{C}_\Q$ can be realised as a 
semigroup crossed product $C^*(\Q/\Z)\rtimes_\alpha \N^*$ by an endomorphic  
action $\alpha$ of the multiplicative semigroup $\N^*$ on the group 
$C^*$-algebra 
$C^*(\Q/\Z)$ \cite{LR2}, and this realisation has provided useful 
insight into the analysis of $\mathcal{C}_\Q$ \cite{L2,N}. Since 
individual elements of $\Q/\Z$ and $\N^*$ involve only finitely many 
primes, $C^*(\Q/\Z)\rtimes_\alpha \N^*$ is the direct limit of 
subalgebras $C^*(G_F)\rtimes_\alpha \N^F$, where $F$ is a finite set 
of primes, $G_F$ is the subgroup of $\Q/\Z$ in which the denominators 
have all prime factors in $F$, and $\N^F$ acts through the embedding 
$(n_p)\mapsto\prod_{p\in F}p^{n_p}$ of $\N^F$ in $\N^*$ (see 
Section~\ref{introcrap}). One can therefore hope to understand the 
Hecke algebra $\mathcal{C}_\Q$ in terms of the finite-prime analogues  
$C^*(G_F)\rtimes_\alpha \N^F$. 

Our goal is to analyse the structure of these finite-prime 
analogues of the Bost-Connes algebra. We started this analysis in 
\cite{LPR}, where we described a composition series for the two-prime 
analogue and identified the subquotients in familiar terms: there is a 
large type~I ideal, a commutative quotient isomorphic to $C(\T^2)$, 
and the intermediate subquotient is isomorphic to a direct sum of 
Bunce-Deddens algebras. Here we describe a composition series for 
$C^*(G_F)\rtimes_\alpha \N^F$. Again there are a large type~I ideal 
and a commutative quotient, and the intermediate subquotients are 
direct sums of simple $C^*$-algebras. We can describe the simple 
summands as ordinary crossed products by actions of $\Z^S$ for 
$S\subset F$. When $\vert S\vert =1$, these actions are odometers and 
the crossed products are Bunce-Deddens algebras; when $\vert S\vert 
>1$, the actions are an apparently new class of higher-rank odometer 
actions, and the crossed products are an apparently new class of 
classifiable AT-algebras.

We begin with a short section  in which we describe the algebras we intend to 
study. In 
\S\ref{finitely_many_primes}, we describe our composition series for 
the semigroup crossed product 
$C^*(G_F)\rtimes_\alpha \N^F$. It has $\vert F\vert+1$ subquotients, 
and all but two of them are direct sums of algebras stably isomorphic 
to ordinary crossed products of the form $C(\U(\Z_{F\setminus 
  S}))\rtimes\Z^S$, where $S\subset F$ and $\U(\Z _{F\setminus S})$ is 
the group of units in the ring $\prod_{p\in F\setminus S}\Z_p$. Our 
main tools are the analysis of invariant ideals in semigroup 
crossed products from \cite{Lar} and some technical lemmas on sums and 
intersections of 
ideals in $C^*$-algebras. We also use the general results of \cite{P} 
to see that the simple summands are classifiable. 

In \S\ref{section_one_prime_on_two}, we show that when $S=\{q\}$ is a 
singleton, $C(\U(\Z_{F\setminus S}))\rtimes\Z^S$ is a direct sum of 
finitely many Bunce-Deddens algebras; as in \cite{LPR}, the number of summands  
depends 
on the orders of $q$ in the finite groups $\prod_{p\not=
  q}\U(\Z/p^{l}\Z)$ for large $l\in \N$. We then consider the case 
where $S=\{q,r\}$. By computing the $K$-theory of  $C(\U(\Z_{F\setminus
  S}))\rtimes\Z^S$, we can see that they are not 
Bunce-Deddens algebras, for example. We expect these summands to be 
even harder to recognise when $\vert S\vert>2$. 

In \S\ref{BCalgebra}, we use techniques like those of 
\S\ref{finitely_many_primes} to identify 
subquotients of the Bost-Connes algebra $C^*(\Q/\Z)\rtimes_\alpha
\N^*$. These include algebras stably isomorphic to 
$C(\U(\Z_{\mathcal{P}\setminus S}))\rtimes \Z^S$ when $S$ is a cofinite 
subset of the set $\mathcal{P}$ of all primes; in this case, though, these 
crossed products  are 
themselves simple, and even though the general theory of \cite{P} no 
longer applies, we can see using results from \cite{BDR} that they are classifiable 
AT-algebras. We finish with a purely 
number theoretic 
Appendix in which we identify the orders of an odd integer in the 
groups $\U(\Z/p^l\Z)$ and their products. As in 
\cite[Theorem~3.1]{LPR}, these are needed when we want to identify the number of 
simple summands in the various subquotients.

\subsection*{Acknowledgements} 
We thank Chris Phillips for pointing out that the AH-algebras appearing as simple  
summands in our structure theorems are in fact AT-algebras.

\section{Preliminaries}\label{introcrap}

We denote by $\N^*$ the semigroup of positive integers under 
multiplication, and by $\N$ the semigroup of nonnegative integers 
under addition. It was shown in \cite[Proposition~2.1]{LR2} that there 
is an action $\alpha$ of $\N^*$ by endomorphisms of $C^*(\Q/\Z)$ such that 
\[
\alpha_n(\delta_r)=\frac 1n \sum_{ns=r}\delta_s\ \text{ for $r\in
  \Q/\Z$ and $n\in\N^*$.}
\]
The corresponding semigroup crossed product  $C^*(\Q/\Z)\rtimes_\alpha 
\N^*$ is isomorphic to the Hecke $C^*$-algebra $\mathcal{C}_\Q$ of 
Bost and Connes \cite[Corollary~2.10]{LR2}. We denote by 
$(i_A,i_{\N^*})$ the canonical covariant representation of 
$(C^*(\Q/\Z),\N^*,\alpha)$ in $C^*(\Q/\Z)\rtimes_\alpha 
\N^*$.

Let $F$ be a set of prime numbers.  The rational numbers of the form 
$k\big(\prod_{p\in F}p^{m_p}\big)^{-1}$ form a subgroup of $\Q$, whose 
image in $\Q/\Z$ we 
denote by $G_F$. The integrated form of the map 
$r\mapsto \delta_r:G_F\to UC^*(\Q/\Z)$ is a homomorphism $i_F$ of $C^*(G_F)$ 
into $C^*(\Q/\Z)$; a standard duality argument shows that $i_F$ is 
injective, so that we can identify $C^*(G_F)$ with the subalgebra 
$i_F(C^*(G_F))$ of $C^*(\Q/\Z)$. When $n$ has all its prime factors in 
$F$, $\alpha_n$ leaves this subalgebra invariant, and hence composing 
$\alpha$ with the map $(m_p)_{p\in F}\mapsto \prod_{p\in F}p^{m_p}$ 
gives an action of $\N^F$ on $C^*(G_F)$, which we also denote by 
$\alpha$. The pair $(i_F,i_{\N^*}\vert_{\N^F})$ is a covariant 
representation of $(C^*(G_F),\N^F,\alpha)$ in $C^*(\Q/\Z)\rtimes_\alpha 
\N^*$. Since $i_F$ is injective, we can deduce from the main theorem 
of \cite{LarR} (or by minor modifications to the argument in \S3 
of \cite{LR2}) that the corresponding homomorphism 
\[
i_F\times i_{\N^*}\vert_{\N^F}:C^*(G_F)\rtimes_\alpha \N^F\to 
C^*(\Q/\Z)\rtimes_\alpha\N^*
\]
is also an injection. We use this injection to identify 
$C^*(G_F)\rtimes_\alpha \N^F$ with a subalgebra of $C^*(\Q/\Z)\rtimes_\alpha\N^*$.

The crossed product $C^*(\Q/\Z)\rtimes_\alpha\N^*$ is spanned by the 
elements of the form $i_A(\delta_r)i_{\N^*}(m)i_{\N^*}(n)^*$ 
\cite[Lemma~3.2]{LR2}. If $F$ contains all the prime factors of 
$m$, $n$ and the denominator of $r$, then 
$i_A(\delta_r)i_{\N^*}(m)i_{\N^*}(n)^*$ lies in 
$C^*(G_F)\rtimes_\alpha \N^F$. Thus   $C^*(\Q/\Z)\rtimes_\alpha\N^*$ 
is the direct limit $\overline{\bigcup_{F}C^*(G_F)\rtimes_\alpha 
  \N^F}$ over increasing finite subsets $F$ of the 
set $\mathcal P$ of prime numbers. 

In the next section, we shall describe a composition series for  
$C^*(G_F)\rtimes_\alpha \N^F$ when $F$ is a finite subset of $\mathcal 
P$, and identify the subquotients in terms 
of ordinary crossed products $C(X_S)\rtimes \Z^S$ associated to subsets 
$S$ of $F$. The underlying space $X_S$ is the group of units 
$\U(\Z_{F\setminus S})$ in the ring $\Z_{F\setminus S}:=\prod_{p\in 
  F\setminus S}\Z_p$; as an additive group, $\Z_{F\setminus S}$ is the 
dual group of $G_{F\setminus S}$. The action of a prime $q\in S$ on  
\[
C(\U(\Z_{F\setminus S}))\subset C(\Z_{F\setminus S})\cong 
C^*(G_{F\setminus S})\subset C^*(\Q/\Z)
\]
induced by $\alpha_q$ is multiplication by $q$ on $\U(\Z_{F\setminus S})$  
(see \cite[Lemma~1.1]{LPR}), which because $q$ is a unit in 
$\Z_{F\setminus S}$ is an automorphism. Thus the action of $\N^S$ on 
$C(\U(\Z_{F\setminus S}))$ extends to an 
action $\sigma$ of $\Z^S$ such that 
\[\sigma_{(m_p)}(f)(x)= f\big(\big(\textstyle{\prod_{p\in S}}p^{m_p}\big)^{-1}x\big)  
\text{ for }(m_p)\in \N^S.
\]

As a matter of notation, 
we shall view a crossed product $A\rtimes_\beta G$ by an action of a group as 
the universal 
$C^*$-algebra generated by a copy of $A$ and a unitary 
representation $i_G:G\to U(A\rtimes_\beta G)$ satisfying the 
covariance relation $\beta_s(a)=i_G(s)ai_G(s)^*$.

\section{Finitely many primes}\label{finitely_many_primes}

The object of this section is to prove the following theorem. For the 
definitions of AT-algebra, real rank zero and stable rank one, see \cite{Ro}  
and the references given there.

\begin{theorem}\label{comp_series}
Let $F$ be a finite set of primes. Then there is a composition series 
$\{I_k\mid 1\leq k \leq \vert F \vert\}$ of ideals in 
$C^*(G_F)\rtimes_\alpha \N^F$ such that 
\begin{enumerate}\smallskip
\item[(a)] $I_1\cong C(\U(\Z_F), \K(l^2(\N^F)))$;
\smallskip
\item[(b)] $I_{k+1}/I_k \cong \bigoplus_{S\subset F,\,\vert S
    \vert=k}\bigl(C(\U(\Z_{F\setminus S}))\rtimes_{\sigma} \Z^S  
  \bigr)\otimes \K(l^2(\N^{F\setminus S}))$;
\smallskip
\item[(c)] $(C^*(G_F)\rtimes\N^F)/I_{\vert F\vert}\cong
  C(\T^F).$
\smallskip
\end{enumerate}
Each $C(\U(\Z_{F\setminus S}))\rtimes_{\sigma} \Z^S$ is a finite 
direct sum  of simple AT-algebras with 
real rank zero and a unique tracial state.
\end{theorem}

The proof of the theorem will occupy the rest of the section. We 
will need some notation and a number of preliminary results.

Under the Fourier transform $C^*(G_F)\cong C(\Z_F)$ the action 
$\alpha$ becomes
$$
\alpha_{(n_p)}(f)(x)=
\begin{cases}f\big(\big(\textstyle{\prod_{p\in
      F}}p^{n_p}\big)^{-1}x)&\text{if 
$x\in \big(\textstyle{\prod_{p\in F}}p^{n_p}\big)\Z_F$}\\
0&\text{otherwise}
\end{cases}
$$
(see \cite[Lemma~1.1]{LPR}). For $S\subset F$, we set $\z_S:=\{a\in \Z_F\mid 
a_p=0\text{ for }p\in S\}$,  and we write 
$\z_p$ for $\z_{\{p\}}$. The next lemma identifies $C_0(\Z_F\setminus 
\z_S)$ as the kind of ideal for which taking crossed products behaves 
well (see \cite{Lar}).

\begin{lemma} For $S\subset F$, 
$C_0(\Z_F\setminus \z_S)$ is an extendibly invariant ideal in 
$C(\Z_F).$
\label{ext_ideals_of_BC_algebra}
\end{lemma}

\begin{proof} It suffices by \cite[Theorem~4.3]{Lar} to show that for 
  each $n\in \N^F$, the 
  endomorphism $x\mapsto \big(\prod_{p\in F}p^{n_p}\big)x$ of $\Z_F$   
  leaves both   $\z_S$ and $\Z_F\setminus \z_S$ 
  invariant. Certainly $\big(\prod_{p\in F}p^{n_p}\big)\z_S$ is 
  contained in $\z_S$. If 
  $x\notin \z_S$, then 
  $x_r\not=0$ for some $r\in S$, $\prod_{p\in F}p^{n_p}x_r\not=0$ for 
  this $r$, and   $\big(\prod_{p\in F}p^{n_p}\big)x\notin   \z_S$.  
\end{proof}

Theorem~1.7 of \cite{Lar} now allows us to identify $C_0(\Z_F\setminus 
\z_S)\rtimes \N^F$ with an ideal $J_S$ in $C(\Z_F)\rtimes_\alpha\N^F$ 
such that $(C(\Z_F)\rtimes_\alpha\N^F)/J_S=C(\z_S)\rtimes \N^F$; we 
write $J_p$ for $J_{\{p\}}$.  

\begin{lemma}$J_S=\sum_{p\in S}J_p$. 
\label{J_S_as_sum}
\end{lemma}

\begin{proof} Since $\z_S=\bigcap_{p\in S}\z_p$, we have  
  $\Z_F\setminus \z_S=\bigcup_{p\in S}\Z_F\setminus \z_p$, and  
$C_0(\Z_F\setminus \z_S)=\sum_{p\in S}C_0(\Z_F\setminus 
\z_p)$. It follows from Lemma 1.3 of \cite{Lar}  that if $I$, $J$ and 
$I+J$ are extendibly invariant ideals in $(A, P)$, then $(I+J)\rtimes 
P=(I\rtimes P)+(J\rtimes P)$. Thus the result follows from 
Lemma~\ref{ext_ideals_of_BC_algebra}.
\end{proof}

For $1\leq k\leq \vert F\vert$, we define 
\begin{equation}
I_k:= \textstyle{\prod_{S\subset F,\,\vert S
    \vert=k}}J_S=\textstyle{\bigcap_{S\subset F,\,\vert S
    \vert=k}}J_S. 
\label{def_of_Ik}
\end{equation}
It follows from \cite[Lemma~1.3]{Lar} that if $I$ and $J$ are 
extendibly invariant ideals in $(A,P)$, then
\[
(I\rtimes P)(J\rtimes P)=(IJ)\rtimes P, 
\] 
and hence  
$I_k=C_0\big(\textstyle{\bigcap_{S\subset F,\,\vert S
    \vert=k}}(\Z_F\setminus \z_S)\big)\rtimes \N^F$. Therefore 
\[
I_1=C_0\big( \textstyle{\bigcap_{p\in F}}(\Z_F\setminus
\z_p)\big)\rtimes \N^F 
=C_0\big(\textstyle{\prod_{p\in F}}(\Zp \setminus
\{0\})\big)\rtimes \N^F; 
\]
since $\Zp \setminus \{0\}$ is homeomorphic to  $\U(\Zp)\times \N$ 
by \cite[Lemma 2.3]{LPR}, part (a) of Theorem~\ref{comp_series} follows 
from an argument similar to the one in the last paragraph of \cite[page 
176]{LPR}. Similarly, we can prove part (c) by following the proof of 
(2.4) of 
\cite{LPR}, because $\big(C^*(G_F)\rtimes_\alpha \N^F\big)/I_{\vert 
  F\vert}=\C\rtimes \N^F$.

To prove part (b) of Theorem~\ref{comp_series}, we need some lemmas. The first contains 
some 
general facts about families of ideals in $C^*$-algebras.

\begin{lemma} 
Suppose that $I_1, \dots ,I_n$ are ideals in a $C^*$-algebra $B$. 
\begin{enumerate} 
\item[(a)] With $F_n=\{1, \dots, n\}$, we have 
\begin{equation} 
{\prod_{S\subset F_n,\,\vert S \vert=k}}\big({\textstyle{\sum_{i\in 
    S}}} I_i\big)={\sum_{R\subset F_n,\,\vert 
    R \vert=n-k+1}}\big(\textstyle{\prod_{j\in R}}I_j\big) \text{ for } 1\leq 
k \leq n. 
\label{prod_sum_is_sum_prod} 
\end{equation} 
\item[(b)] 
Suppose that $K$ is an ideal such that $I_iI_j\subset K$ for all 
$i, j$. Then $(\sum_{i=1}^nI_i)/K$ is naturally isomorphic to 
$\bigoplus_{i=1}^n(I_i/I_i\cap K).$ 
\smallskip 
\end{enumerate} 
\label{lemma_prod_sum_is_sum_prod} 
\end{lemma} 
 
\begin{proof} 
We prove (a) by induction on $n$. The statement is trivial for 
$n=1,2$. Suppose it holds for $n-1$. When $k=1$, both sides of 
(\ref{prod_sum_is_sum_prod}) are $\prod_{i=1}^nI_i$, so we assume 
$k\geq 2$. Writing the left-hand side $\LHS$ of (\ref{prod_sum_is_sum_prod}) 
as $(\prod_{n\in S})(\prod_{n\notin S})$ and applying the inductive 
hypothesis to $F_{n-1}$ shows that  
\begin{equation} 
\LHS=\Big({\prod_{\vert S\vert=k, n\in S}}\big(I_n+{\textstyle{\sum_{i\in 
  S\setminus\{n\}}}}I_i\big)\Big) \Big({\sum_{R\subset F_{n-1} 
,\,\vert  R \vert=n-k}}\big(\textstyle{\prod_{j\in R}}I_j\big) \Big). 
\label{LHS_of_step_n} 
\end{equation} 
Because $I_n$ is an ideal and $I_n^2=I_n$, the first 
term of (\ref{LHS_of_step_n}) simplifies to give  
$$ 
\LHS=\Big(I_n + {\prod_{S'\subset F_{n-1},\,\vert S' 
    \vert=k-1}} \big({\textstyle{\sum_{i\in 
  S'}}}I_i\big)\Big)\Big(\sum_{R\subset F_{n-1} 
  ,\,\vert R\vert=n-k}\big(\textstyle{\prod_{j\in R}}I_j\big) \Big). 
$$ 
We can use the inductive hypothesis on $F_{n-1}$ with $k$ 
replaced by $k-1$ to get 
\begin{equation} 
\LHS=\Big(I_n+ {\sum_{R'\subset F_{n-1},\,\vert 
    R'\vert=n-k+1}}\big({\textstyle{\prod_{j\in R'}}}I_j\big) \Big)  
\Big({\sum_{R\subset F_{n-1},\,\vert 
    R\vert=n-k}}\big({\textstyle{\prod_{j\in R}}}I_j\big) \Big), 
\label{big_mess_one} 
\end{equation} 
which is contained in 
\begin{equation} 
{\sum_{R\subset F_{n-1},\,\vert 
    R\vert=n-k}}\big({\textstyle{\prod_{j\in R\cup\{n\}}}}I_j) +  
{\sum_{R'\subset F_{n-1},\,\vert 
    R'\vert=n-k+1}}\big({\textstyle{\prod_{j\in R'}}}I_j\big). 
\label{big_mess_two} 
\end{equation} 
Since (\ref{big_mess_two}) 
is the same as the right-hand side $\RHS$ of 
(\ref{prod_sum_is_sum_prod}), $\LHS\subset \RHS$.  
On the other hand, every element of 
every $\prod_{j\in R'}I_j$ arises in (\ref{big_mess_one}) because we 
can pick $R\subset R'$, so $\RHS\subset \LHS$.

To prove (b), note that the map $\phi_i:a+I_i\cap K\mapsto a+K$ is an 
injection of $I_i/(I_i\cap K)$ into $\big(\sum_{i=1}^nI_i\big)/K$, and 
$$
\phi_i(a+I_i\cap K)\phi_j(b+I_j\cap K)=ab+K=0\text{ for }i\neq j 
$$ 
because $ab\in I_iI_j\subset K$. So the $\phi_j$ combine to give an 
injection $\phi$ of $\bigoplus(I_i/I_i\cap K)$ into 
$(\sum_{i=1}^nI_i)/K$, which is clearly surjective. 
\end{proof}

\begin{lemma}\label{k_subquotient_as_scp} The ideals $I_k$ of 
  $C^*(G_F)\rtimes_\alpha \N^F$ defined in 
  \textnormal{(\ref{def_of_Ik})} satisfy 
$$
I_{k+1}/I_k=\bigoplus_{S\subset F,\,\vert S \vert=k}\big(\textstyle{\bigcap_{p\notin 
  S}}J_{S\cup \{p\}}\big)/J_S. 
$$
\end{lemma}

\begin{proof}
Lemma~\ref{lemma_prod_sum_is_sum_prod} (a) gives 
$I_{k+1}=\textstyle{\sum_{R\subset F,\,\vert R\vert=n-k}}(\prod_{p\in 
  R}J_p)$. The product of any 
two ideals  $\prod_{p\in R}J_p$ with $\vert R \vert = n-k$ has at 
least $n-k+1$ factors $J_p$, and hence is contained in 
$I_k=\textstyle{\sum_{R\subset F,\,\vert R\vert=n-k+1}}(\prod_{p\in  
  R}J_{p})$. Thus part (b) of Lemma~\ref{lemma_prod_sum_is_sum_prod} 
gives 
\begin{equation} 
I_{k+1}/I_k= \bigoplus_{R\subset F,\,\vert R\vert=n-k}\frac{\prod_{p\in 
    R}J_p}{I_k\cap\big(\prod_{p\in R}J_p\big)}.  
\label{subquotient_as_direct_sum_1} 
\end{equation} 
Now  
\[I_k\cap(\textstyle{\prod_{p\in R}}J_p)=\sum_{\vert T \vert 
  =n-k+1}(\textstyle{\prod_{q\in T}}J_q)(\textstyle{\prod_{p\in R}}J_p); 
\]  each of these 
summands has at least one factor $J_q$ for $q\notin R$, and is then 
contained in $J_q(\prod_{p\in R}J_p)$. Using $I\cap J=IJ$ again gives 
\[
I_k\cap\big(\textstyle{\prod_{p\in R}}J_p\big)=\textstyle{\sum_{q\notin 
    R}}J_q\big(\prod_{p\in R}J_p\big)=\big(\textstyle{\sum_{q\notin 
    R}}J_q 
\big)\big(\prod_{p\in R}J_p \big), 
\]
and using the isomorphism $(I+J)/I=J/(I\cap J)$ and 
Lemma~\ref{J_S_as_sum} gives 
\[
\frac{\textstyle{\prod_{p\in R}}J_p}{I_k\cap(\prod_{p\in R}J_p)}= 
\frac{J_{F\setminus R}+ \big(\prod_{p\in R}J_p \big)}{J_{F\setminus 
    R}}. 
\]
Finally we observe that 
\[
J_{F\setminus R}+ \big(\textstyle{\prod_{p\in R}}J_p 
\big)=\textstyle{\prod_{p\in R}}\big(J_{F\setminus R}+J_p 
\big)=\textstyle{\prod_{p\in R}}J_{(F\setminus   R)\cup\{p\}} 
\] and write $S=F\setminus R$ to deduce the result. 
\end{proof}

\begin{lemma}\label{identify_small_subquotient_as_scp} 
The ideals $J_S$ in $C^*(G_F)\rtimes_\alpha \N^F$ satisfy 
$$\big(\textstyle{\bigcap_{p\in F\setminus S}}J_{S\cup\{p\}}\big)/J_S\cong 
\big(C(\U(\Z_{F\setminus S}))\rtimes_\sigma \Z^S\big) 
\otimes \K(l^2(\N^{F\setminus S})).
$$
\end{lemma}

\begin{proof} We first realise $\big(\bigcap_{p\in F\setminus 
    S}J_{S\cup\{p\}}\big)/J_S$ as a semigroup crossed product: 
\begin{align*}
\textstyle{\bigcap_{p\in F\setminus S}}J_{S\cup\{p\}} 
&=C_0\big(\textstyle{\bigcap_{p\in F\setminus S}}(\Z_F\setminus \z_{S\cup 
  \{p\}}\big)\big)\rtimes \N^F\\ 
&= 
C_0\big(\Z_F\setminus\big(\textstyle{\bigcup_{p\in F\setminus S}}\z_{S\cup
  \{p\}}\big)\big)\rtimes \N^F. 
\end{align*}
Thus  
\begin{align*}
\big({\textstyle{\bigcap_{p\in F\setminus S}}J_{S\cup\{p\}}}\big)/{J_S} 
&=C_0\big(\z_S\setminus\big(\textstyle{\bigcup_{p\in F\setminus S}}\z_{S\cup 
  \{p\}}\big)\big)\rtimes\N^F\\ 
&=C_0\big(\big(\textstyle{\prod_{p\in F\setminus 
    S}}\Z_p\setminus\{0\}\big)\times 
\big(\textstyle{\prod_{p\in S}}\{0\}\big)\big)\rtimes\N^F. 
\end{align*}
The arguments of Corollary~2.4 and Lemma~2.5 of \cite{LPR} show that 
this last crossed product is isomorphic to 
$\big(C(\U(\Z_{F\setminus S}))\rtimes_\sigma \Z^S\big) 
\otimes \K(l^2(\N^{F\setminus S}))$.
\end{proof}

Part (b) of Theorem~\ref{comp_series} follows immediately from 
Lemmas~\ref{k_subquotient_as_scp} and 
\ref{identify_small_subquotient_as_scp}.

To finish the proof of Theorem~\ref{comp_series}, it remains to prove 
the statements about the structure of $C(\U(\Z_{F\setminus 
  S}))\rtimes_\sigma \Z^S$. Corollary~\ref{indexform} implies that 
$H:=\overline{\Z^S}$ has 
finite index in $\U(\Z_{F\setminus 
  S})$. The argument at the end of the proof of 
\cite[Theorem~3.1]{LPR} shows that 
$C(\U(\Z_{F\setminus S}))\rtimes_\sigma \Z^S$ is a finite direct sum 
of algebras isomorphic to $C(H)\rtimes_\sigma \Z^S$, which is simple 
because $\Z^S$ acts minimally and freely on $H$. Since $H$ is an open 
and closed subset of $\U(\Z_{F\setminus S})$, it is totally 
disconnected, and it follows from \cite[Theorem~6.11]{P} that  
$C(H)\rtimes_\sigma
 \Z^S$ has real rank zero and stable rank one.

The space $ \U(\Z_{F\setminus S})$ is the inverse limit of 
the finite groups $\U\big(\Z/\big(\prod_{p\in F\setminus S} 
p^{l_p}\big)\Z\big)$ over $ l=(l_p)\in 
\N^{F\setminus S}$. The diagonally embedded copy of $\N$ is 
cofinal in $\N^{F\setminus S}$, and hence 
\begin{equation}
\U(\Z_{F\setminus S})=\varprojlim
\U\big(\Z/\big(\textstyle{\prod_{p\in 
    F\setminus S}} p^{n}\big)\Z\big).
\label{units_as_inverse_limit}
 \end{equation}
Let $\pi_n$ denote the canonical surjection of 
$\U(\Z_{F\setminus S})$ onto $\U\big(\Z/\big(\textstyle{\prod_{p\in
    F\setminus S}} p^{n}\big)\Z\big)$.

\begin{lemma} 
Let $H_n:= \pi_n(H)\subset \U\big(\Z/\big(\textstyle{\prod_{p\in 
    F\setminus S}}  p^{n}\big)\Z\big)$ and let $\Z^S$ act on $H_n$ via 
the embedding $(n_q)\mapsto 
\prod_{q\in S}q^{n_q}$  of $\Z^S$ in $\Z$. Then there are 
$C^*$-subalgebras $A_n$ of  $C(H)\rtimes_\sigma\Z^S$ such that 
$A_n\cong C(H_n)\rtimes \Z^S$ and $C(H)\rtimes_\sigma 
\Z^S=\overline{\bigcup A_n}$. 
\label{simple_bit_as_union_of_Zd_cp} 
\end{lemma}

\begin{proof} The homomorphism $\pi_n$ induces an injection $\pi_n^*$ 
  of $C(H_n)$ into $C(H)$, and then 
$C(H)= 
\overline{\bigcup_{n\in \N}\pi_n^*(C(H_n))}$. 
On $\Z^S\subset H$, $\pi_n$ is reduction 
modulo $\prod_{p\in  F\setminus S}p^{n}$, so $\pi_n^*$ 
converts the action $\sigma$ into the canonical action of $\Z^S$ by 
multiplication on $H_n$. Thus $\pi_n^*$ induces a homomorphism 
$\pi_n^*\rtimes \operatorname{id}$ of $C(H_n)\rtimes \Z^S$ into 
$C(H)\rtimes_\sigma \Z^S$. The 
homomorphism 
$\pi_n^*$ is faithful on $C(H_n)$ and intertwines the dual actions, 
and hence a standard argument shows that $\pi_n^*\rtimes 
\operatorname{id}$  
is faithful on $C(H_n)\rtimes \Z^S$  (see, for example, \cite[Lemma~4.2]{LPR}). 
Since $\bigcup_n \pi_n^*(C(H_n))$ is dense in 
$C(H)$, we therefore have 
\[
C(H)\rtimes_\sigma \Z^S 
=\overline{\textstyle{\bigcup_{n\in \N}}\pi_n^*\rtimes 
  \id (C(H_n)\rtimes \Z^S)}, 
\]
as claimed.
\end{proof}

We can identify the  subalgebras $A_n$ explicitly.

\begin{prop}\label{Smithnf} 
Let $F$ be a finite quotient of $\Z^k$. Then 
  $C(F)\rtimes \Z^k$ is isomorphic to $C(\T^k, M_{\vert F\vert}(\C)).$ 
\label{unwinding_Zn_cp}
\end{prop}

\begin{proof} 
Let $H$ be the subgroup of $\Z^k$ with $F= \Z^k/H$. Then $H$ is itself 
a free abelian group of rank $k$, and hence has the form $A\Z^k$ for 
some $A\in M_k(\Z)\cap GL_k(\Q)$. The matrix $A$ has a Smith normal 
form: there are matrices $P, Q\in GL_k(\Z)$ such that $B:=P^{-1}AQ^{-1}$ 
is diagonal \cite[\S3.22]{MM}. Then $H=A\Z^k=PBQ\Z^k=PB\Z^k\cong 
B\Z^k=b_{11}\Z\oplus \dots \oplus b_{kk}\Z$. In other words, multiplying by 
$P^{-1}$ gives an automorphism of $\Z^k$ which carries $H$ into 
$\bigoplus b_{ii}\Z$. Thus 
\[
C(F)\rtimes \Z^k 
\cong C\big(\textstyle{\prod_{i=1}^k}(\Z/b_{ii}\Z)\big)\rtimes \Z^k 
\cong \textstyle{\bigotimes_{i=1}^k}\big(C(\Z/b_{ii}\Z)\rtimes_\tau \Z\big), 
\] 
where $\tau$ is the action of $\Z$ by translation. 

By \cite[Corollary 2.5]{OP}, $C(\Z/b\Z)\rtimes_\tau \Z$ is isomorphic 
to the induced algebra 
\[
\Ind_{(b\Z)^\perp}^\T(C(\Z/b\Z)\rtimes_\tau 
(\Z/b\Z),\;\widehat\tau), 
\]
which is described in terms of a generator $\beta$ of the dual 
action $\widehat{\tau}$ as the mapping torus 
\begin{equation} 
MT(\beta)=\{f:[0, 1]\to C(\Z/b\Z)\rtimes_\tau (\Z/b\Z) \mid 
f(1)=\beta(f(0))\}. 
\label{mapping_torus} 
\end{equation}
Since $C(\Z/b\Z)\rtimes_\tau (\Z/b\Z)\cong B(l^2(\Z/b\Z))=M_{\vert b\vert 
  }(\C)$, the automorphism $\beta$ is inner, and there is a continuous 
path $\beta_t$ in $\Aut M_{\vert b\vert}(\C)$ such that $\beta_0=\id$ 
and $\beta_1=\beta$. Now $\phi(f)(t)=\beta_t^{-1}(f(t))$ defines an 
isomorphism $\phi$ of (\ref{mapping_torus}) onto $C(\T, M_{\vert 
  b\vert}(\C))$. We therefore have 
\[
C(F)\rtimes\Z^k \cong \textstyle{\bigotimes_{i=1}^k}C(\T, M_{\vert 
  b_{ii}\vert}(\C))\cong C(\T^k, M_{\prod_i \vert b_{ii}\vert}(\C)), 
\] 
and the result follows on observing that $\prod_i\vert b_{ii}\vert 
=\vert \det B\vert=\vert \det A\vert=\vert F\vert.$ 
\end{proof}

It follows from Proposition~\ref{Smithnf} and the decomposition 
$C(H)\rtimes_\sigma \Z^S=\overline{\bigcup A_n}$ that 
$C(H)\rtimes_\sigma \Z^S$ 
is an AH-algebra\footnote{To see that an inductive limit 
  $\overline{\bigcup A_n}$  is  
an  AH-algebra, it suffices to show that each $A_n$ is a corner in a matrix 
  algebra $M_N(C(X))$, or, equivalently, that $A_n$ is a homogeneous 
  algebra with vanishing Dixmier-Douady class. Since the 
  Dixmier-Douady class $\delta(A)$ of an $m$-homogeneous algebra 
  satisfies $m\delta(A)=0$, and $H^3(\T^k,\Z)$ has no torsion, it 
  suffices to prove that each $A_n$ is a homogeneous algebra with 
  spectrum $\T^k$. In our situation we could prove this in several 
  ways. However, Proposition~\ref{unwinding_Zn_cp} makes the stronger 
  statement that $A_n$ is isomorphic to $M_m(C(\T^k)).$}. 
The K-theory of $C(H_n)\rtimes \Z^S$ is torsion-free  and  
this property is preserved under inductive limits, so  
$C(H)\rtimes_\sigma \Z^S$ has torsion-free K-theory. Since  
$C(H)\rtimes_\sigma \Z^S$ is a simple AH-algebra with real rank 
zero and no dimension growth, it follows from results of Elliott and Lin 
that it is an AT-algebra (see \cite[Lemma~7.5]{P2}).

We also use the decomposition $C(H)\rtimes_\sigma 
\Z^S=\overline{\bigcup A_n}$ to prove that $C(H)\rtimes_\sigma \Z^S$ 
has a unique tracial state. Let $\mu$ denote the Haar measure on 
$H\subset \U(\Z_{F\setminus 
  S})$. The  
action $\sigma$ permutes the cylinder sets $\{\pi_n^{-1}(m)\mid m\in 
H_n\}$, so every invariant probability measure agrees with $\mu$ on 
cylinder sets. Since the characteristic functions of such sets span a 
dense subspace of $C(H)$, it 
follows that $\mu$ is the only invariant probability measure, and 
$C(H)\rtimes \Z^S$ has a unique tracial 
state by \cite[Corollary~VIII.3.8]{D}. 

This completes the proof of 
Theorem~\ref{comp_series}.

\section{The structure of $C(\U(\Z_{F\setminus S}))\rtimes_\sigma 
  \Z^S$}\label{section_one_prime_on_two} 

\subsection{When $S$ contains just one prime}\label{oneonF} 
We consider $C(\U(\Z_{F\setminus S}))\rtimes_\sigma 
  \Z^S$ when $S=\{q\}$. To simplify the notation, we relabel 
$F\setminus \{q\}$ as $F$. The following result generalises \cite[Theorem 3.1]{LPR} in 
two directions: to sets $F$ with $\vert F \vert>1$ and to sets $F$ 
containing the even prime $2$. If ${\bf l}=(l_p)\in \N^F$ is a multi-index, 
we write $\order_{\bf l}(q)$ for the order of $q$ in $\prod_{p\in 
  F}\U(\Z/p^{l_p}\Z)$. 

\begin{theorem}\label{one_prime_on_n} 
Suppose $F$ is a finite set of primes and $q$ is a prime which does 
not belong to $F$. Then there are a multi-index ${\bf K}=(K_p)\in 
\N^F$ and $d\in \N$ such that 
\begin{equation}\label{ordermodF}
\order_{{\bf K}+{\bf l}}(q)=d\big(\textstyle{\prod_{p\in F}p^{l_p}}\big)
\ \text{for every ${\bf l}\in \N^F$,}
\end{equation}
and $C(\U(\Z_{F}))\rtimes_{\sigma}\Z$ 
is the direct sum of $\big(\prod_{p\in F}(p-1){p}^{K_p-1}\big)/d$
copies of a 
Bunce-Deddens algebra with supernatural number $d\bigl({\prod_{p\in
    F}}p^{\infty}\bigr)$. 
\end{theorem}

The existence of ${\bf K}$ and $d$ satisfying
(\ref{ordermodF}) is established in Proposition~\ref{nothyfacts}. We saw in \S\ref{finitely_many_primes} that
$C(\U(\Z_{F}))\rtimes_{\sigma}\Z$ is the direct sum of copies of the
simple algebra $C(H)\rtimes_\sigma \Z$, where $H$ is the closure of
$q^\Z$ in $\U(\Z_F)$. It remains to prove that $C(H)\rtimes_\sigma 
\Z$ is a Bunce-Deddens algebra and to calculate the index $\vert
\U(\Z_F):H\vert$, which is the number of simple direct summands.

Let $\{n_k\}$ be integers with $n_k\geq 2$,  and
let $X_k=\{0, 1, \dots, n_k -1\}$. Addition by $1$ with carry over is
a homeomorphism of the totally disconnected space $X:=\prod_{k\geq
  0}X_k$ called an odometer action, and the resulting crossed product 
$C(X)\rtimes_\tau \Z$ is a
\emph{Bunce-Deddens algebra with supernatural number} 
${\bf n}:=\prod_{k\geq 0}n_k$ (see \cite[Chapter~VIII.4]{D}).

Our claim that $C(H)\rtimes_\sigma \Z$ is a Bunce-Deddens algebra 
will follow from the next proposition, which generalises \cite[Proposition
3.6]{LPR}.

\begin{prop}
Suppose $\{G_l\mid \in\N\}$ are finite groups and $G=\varprojlim(G_l, \pi_l)$. 
Fix $g\in G$ and let $L$ denote 
the closed subgroup of $G$ generated by $g$. Consider the action 
$\sigma:\Z \to \Aut C(G)$ such that
$\sigma_n(f)(x)=f(g^{-n}x)$. 
Let $\order_l(g)$ denote the order of $\pi_l(g)$ in $G_l$, and let
\begin{equation}
d_l:=\begin{cases} \order_0(g)&\text{ if }l=0\\
\order_l(g)/\order_{l-1}(g)&\text{ if }l\geq 1. \end{cases}
\label{order_of_r}
\end{equation}
Then  $C(L)\rtimes_\sigma \Z$ is 
a Bunce-Deddens algebra with supernatural number $\prod_{l\geq 0}d_l.$
\label{cp_as_BunceDeddens}
\end{prop}

\begin{proof} 
Let $X:=\prod_{l\geq 0}\{0, 1, \dots, d_l-1\}$. 
The argument in the proof of \cite[Proposition~3.6]{LPR} shows that 
the continuous maps $h_l:X\to G_l$ 
given by 
\begin{equation}
h_l(\{a_n\})=\pi_l(g^{a_0+a_1d_0+\dots +a_ld_0d_1\dots d_{l-1}})
\end{equation}
combine to give an equivariant homeomorphism $h:X\to L$ which induces
the required isomorphism. 
\end{proof}

Our subgroup $H$ of $\U(\Z_F)$ is the inverse limit $\varprojlim \pi_l(H)$, where $\pi_l:\U(\Z_F)\to \U(\Z/\big(\prod_{p
\in F}p^{K_p + l}\big)\Z)$ is the canonical surjection. Then 
Proposition~\ref{cp_as_BunceDeddens} and (\ref{ordermodF}) imply that
$C(H)\rtimes_\sigma \Z$ is a Bunce-Deddens algebra with supernatural
number $d(\prod_{p\in F}p)^\infty$ for $d=\order_{{\bf K}}(q)$. By 
Corollary~\ref{indexform}, we have that
\begin{equation}
\vert \U(\Z_F):H \vert=\big(\textstyle{\prod_{p\in F}}(p-1){p}^{K_p-1}\big)/d,
\label{index_of_H_in_U(ZF)}
\end{equation}
which finishes the proof of Theorem~\ref{one_prime_on_n}.

\subsection{When $S$ consists of two primes}
We now analyse $C(\U(\Z_{F\setminus S}))\rtimes \Z^S$
when $S=\{q, r\}$. For simplicity, we consider
only the case $F=\{p, q, r\}$, so that we are
interested in the crossed product $C(\mathcal{U}(\Zp))\rtimes_\sigma
\Z^2$, where 
$$
\sigma_{m, n}(f)(x)=f(q^{-m}r^{-n}x).
$$

\begin{theorem}\label{two_primes_on_one}
The $C^*$-algebra $C(\mathcal{U}(\Zp))\rtimes_\sigma \Z^2$ is a finite
direct sum of copies of a simple AT-algebra $A$ which has real rank
zero, a unique tracial state and $K$-theory
satisfying two short exact sequences:
\begin{equation}
\begin{split}
&0\longrightarrow \Z[p^{-1}]\longrightarrow K_0(A) 
{\longrightarrow}\Z \longrightarrow 0\\
&0\longrightarrow \Z\longrightarrow K_1(A) 
{\longrightarrow}\Z[p^{-1}] \longrightarrow
0.\label{K1_group}
\end{split}
\end{equation}
\end{theorem}

Everything except the assertion about $K$-theory was proved in
Theorem~\ref{comp_series}; the simple $C^*$-algebra $A$ is
$C(H)\rtimes_\sigma \Z^2$, where $H$ is the closure of $q^{\Z}r^{\Z}$
in $\U(\Zp)$. We aim to analyse $C(H)\rtimes_\sigma \Z^2$ by writing it as an
iterated crossed product $(C(H)\rtimes_{\sigma^q} \Z) \rtimes
\Z$. The inside crossed product is not simple unless $q^{\Z}$ is dense
in $H$, and it is helpful to reduce to this case using the following
lemma.

\begin{lemma}\label{Z2_cp_as_ZZ_cp} Let $H_q$ denote the closure of $q^\Z$ in
  $\U(\Zp)$. Then $H_q$ has finite index $I(q)$ in $H$, and hence is
  an open and closed subset of $H$. The inclusion
  of $C(H_q)$ in $C(H)$ induces an isomorphism of
  $C(H_q)\rtimes_{\sigma}(\Z\times I(q) \Z)$ onto the corner
  $\chi_{H_q}(C(H)\rtimes_\sigma \Z^2)\chi_{H_q}$.
\end{lemma}

\begin{proof}
Corollary~\ref{indexform} implies that $H_q$ has finite index in
$\U(\Z_p)$ , so it certainly has finite index in $H$.
The inclusion of $C(H_q)$ in $C(H)$  and the map
\[
(m, I(q)n)\mapsto \chi_{H_q}i_{\Z^2}(m, I(q)n)\chi_{H_q}
\]
form  a
covariant representation of $(C(H_q), \Z\times I(q) \Z, \sigma)$ in
$\chi_{H_q}(C(H)\rtimes_\sigma \Z^2)\chi_{H_q}$, and hence give a
homomorphism 
$$
\phi:C(H_q)\rtimes_{\sigma}(\Z\times I(q) \Z) \to
\chi_{H_q}(C(H)\rtimes_\sigma \Z^2)\chi_{H_q}. 
$$
We can identify $(\Z\times I(q)\Z)^\wedge$ with $\T^2/(\Z\times
I(q)\Z)^\perp=\T^2/(1\times C_{I(q)})$, where $C_n$ denotes the group
of $n$th roots of unity, and then $\phi$ carries the
dual action $\hat{\sigma}_{[w,z]}$ into $\hat{\sigma}_{w,z}$; now a
standard argument implies that $\phi$ is injective (or we could apply
\cite[Corollary 4.3]{R}, for example). We have 
$$
\chi_{H_q}(f i_{\Z^2}(m, n))\chi_{H_q}=(f\chi_{H_q})i_{\Z^2}(m,
n)\chi_{H_q}= i_{\Z^2}(m, n)\sigma_{m,n}^{-1}(f\chi_{H_q})\chi_{H_q}.
$$    
Since the support of $\sigma_{m,n}^{-1}(f\chi_{H_q})$ is contained in
$q^{-m}r^{-n}H_q=r^{-n}H_q$, we have 
$$
\sigma_{m,n}^{-1}(f\chi_{H_q})\chi_{H_q}
=\begin{cases}\sigma_{m,n}^{-1}(f\chi_{H_q})
  &\text{if $I(q)$ divides $n$}\\ 
0&\text{otherwise},\end{cases}
$$ 
and
\begin{align*}
\chi_{H_q}(f i_{\Z^2}(m,
n))\chi_{H_q}&=\begin{cases}i_{\Z^2}(m,n)\sigma_{m,n}^{-1}(f\chi_{H_q}) 
  &\text{if $I(q)$ divides $n$}\\ 
0&\text{otherwise}\end{cases}\\
&=\begin{cases}\phi((f\chi_{H_q})i_{\Z^2}(m,n))
  &\text{if $I(q)$ divides $n$}\\ 
0&\text{otherwise}.\end{cases}
\end{align*}
Thus every $\chi_{H_q}(f i_{\Z^2}(m,n))\chi_{H_q}$ is in the range of
$\phi$, and $\phi$ is surjective.
\end{proof}

\begin{cor} \label{strong_Morita_equiv} 
Define $\gamma:\Z\to \Aut (C(H_q)\rtimes_{\sigma^q}
\Z)$ by 
\begin{equation}
\gamma_m(f i_{\Z}(n))=\sigma_{I(q)m}^r(f)i_{\Z}(n).
\label{formula_for_action_by_I(q)}
\end{equation}  
Then $(C(H_q)\rtimes_{\sigma^q} \Z)\rtimes_\gamma \Z$ is isomorphic to
a full corner in $C(H)\rtimes_{\sigma} \Z^2$.
\end{cor}

\begin{proof}
Theorem 4.1 of \cite{PR} gives
 $C(H_q)\rtimes_\sigma
(\Z\times I(q)\Z)\cong (C(H_q)\rtimes_{\sigma^q}\Z)\rtimes I(q)\Z$, so
the result follows from Lemma~\ref{Z2_cp_as_ZZ_cp} on replacing
$I(q)\Z$ by the isomorphic group $\Z$. 
\end{proof}

The analysis in \S\ref{oneonF} 
shows that $C(H_q)\rtimes_{\sigma^q}\Z$
is a Bunce-Deddens algebra. The $K$-theory of Bunce-Deddens algebras
is well-known. To state the version we need, recall that if ${\bf
  n}=(n_k)_{k\geq0}$ is a sequence with $n_k\geq 2$, then $\Z[{\bf
  n}^{-1}]$ denotes 
the set of rational numbers with denominator $N_k=\prod_{i=0}^kn_i$
for some $k\geq 0$.

\begin{prop}\label{K_theory_of_BD}
Suppose ${\bf n}=(n_k)_{k\geq0}$, $X_k=\{0,\dots,n_k-1\}$, $X=\prod
X_k$ and $\tau:\Z\to\Aut C(X)$ is the associated odometer. Then there
are isomorphisms $\phi_0:K_0(C(X)\rtimes_\tau\Z)\to \Z[{\bf n}^{-1}]$ such
that $\phi_0([\chi_{J(a_0,\dots,a_k)}])=N_k^{-1}$ for each cylinder set
$J(a_0,\dots,a_k)$, and $\phi_1:K_1(C(X)\rtimes_\tau\Z)\to \Z$ such that
$\phi_1(i_{\Z}(1))=1$.   
\end{prop}

\begin{proof}
Because $K_1(C(X))=0$, the Pimsner-Voiculescu sequence for the system 
$(C(X), \Z, \tau)$ reduces to
$$
0\longrightarrow K_1(C(X)\rtimes_\tau
\Z)\overset{\delta}{\longrightarrow} K_0(C(X)) \overset{\id
  -{\tau}_*}{\longrightarrow}
K_0(C(X))\overset{{\id}_*}{\longrightarrow} K_0(C(X)\rtimes_\tau
\Z)\longrightarrow 0.
$$
Now let $C_k=\{J(a_0,\dots,a_k)\}$ be the set of cylinder sets of
length $k+1$, and note that $C(X)=\overline{\bigcup_{k=1}A_k}$, where
$A_k=\lsp\{ \chi_J\mid J\in C_k\}$. Each $\chi_J$ for $J\in C_k$ is the
sum of $n_{k+1}$ basis elements of $A_{k+1}$, so the maps
$[\chi_{J(a_0,\dots,a_k)}]\mapsto {N_k}^{-1}$ of $A_k$ into $\mathbb{R}$
combine to give a homomorphism $\psi_0$ of $K_0(C(X))=\varinjlim{K_0(A_k)}$
into $\mathbb{R}$ with range $\Z[{\bf n}^{-1}]$.  Since the generating 
automorphism
$\tau=\tau_1$ permutes $C_k$, the kernel of  $\psi_0$
is the range of $\id  -{\tau}_*$, and hence $\psi_0$ induces the
required isomorphism $\phi_0$ of $K_0(C(X)\rtimes_\phi \Z)$ 
onto $\Z[{\bf n}^{-1}]$. To verify the statement about $K_1$, recall
that $\delta$ is the coboundary map for the Toeplitz extension of
$C(X)\rtimes_\tau \Z$ (see \cite[\S 2]{PV}), and compute the index of
$[i_\Z(1)]$ in 
$K_0(C(X)\otimes \K)\cong K_0(C(X))$.
\end{proof}

\begin{proof}[Proof of Theorem~\ref{two_primes_on_one}.]
We saw in the proof of Proposition~\ref{cp_as_BunceDeddens} and in the 
paragraph following it 
that the homeomorphism $h$ of $\prod_{k\geq 0}X_k$ onto the subgroup
$H_q$ of $\U(\Zp)$ satisfies
$$
\pi_k(h(\{a_n\}))=\pi_k(q^{a_0 + a_1\order_p(q)+\cdots +
a_k\order_p(q)p^{k-1}})\text{ for }k\geq 0,
$$
and hence carries $J(a_0,\dots,a_k)$ onto $Z(q^{a_0 +
  a_1\order_p(q)+\cdots + a_k\order_p(q)p^{k-1}})$, where 
$$
Z_k(n)=\{x\in \U(\Zp)\mid \pi_k(x)=\pi_k(n)\}.
$$ 
So we deduce from Proposition~\ref{K_theory_of_BD} that there is an
isomorphism $\phi_0$ of $K_0(C(H_q)\rtimes_\sigma \Z)$ onto $\frac
1{\order_p(q)}\Z[p^{-1}]$ such that
$$
\phi_0([\chi_{Z_k(m)}])=\frac 1{\order_p(q)}\frac 1{p^{k}}
$$
for every integer $m$ which lies in $H_q$.

Multiplying by the unit $r^{-I(q)l}$ carries
$Z_k(m)$ into $Z_k(r^{-I(q)l}m)$, and hence $\phi_0 \circ
(\gamma_l)_*=\phi_0$. Thus $(\gamma_l)_*$ is the identity on
$K_0(C(H_q)\rtimes_\sigma \Z)$. It is also the identity on
$K_1(C(H_q)\rtimes_\sigma \Z)$, and hence the Pimsner-Voiculescu
sequence for $((C(H_q)\rtimes_\sigma \Z), \Z, \gamma)$ collapses to
the two short exact sequences
\begin{align}
&0\longrightarrow \frac 1{\order_p(q)}\Z[p^{-1}]\longrightarrow
K_0(C(H_q)\rtimes \Z^2)  
{\longrightarrow}\Z \longrightarrow 0\notag \\
&0\longrightarrow \Z\longrightarrow K_1(C(H_q)\rtimes \Z^2) 
{\longrightarrow}\frac 1{\order_p(q)}\Z[p^{-1}] \longrightarrow
0.  \notag  
\end{align}
From this and Corollary~\ref{strong_Morita_equiv} we can deduce
(\ref{K1_group}); since the isomorphism induced by
Corollary~\ref{strong_Morita_equiv}  scales the
class of $[1]$, we have removed the factor ${\order_p(q)}^{-1}$  by a
further 
scaling to ensure that the final statement does not depend on the order
of factors in our decomposition.
\end{proof}

\begin{remark}
The number of simple summands in Theorem~\ref{two_primes_on_one} is 
$\vert\mathcal{U}(\Zp): H\vert$, and we can compute this using
\cite[Lemma 3.7]{LPR}. For example, if $p$ is odd and $l$ is large, we have 
from (\ref{order_of_m_mod_pl}) that  
\begin{align}
\vert \pi_l(H) \vert
&=[\order_{p^l}(q), \order_{p^l}(r)]=[p^{l-L_p(q)}\order_p(q),
q^{l-L_p(r)}\order_p(r)]\notag \\
&=p^{l-\text{min}(L_p(q),L_p(r))}[\order_p(q),\order_p(r)];\notag
\end{align}
thus we deduce 
$$
\vert\mathcal{U}(\Zp): H\vert=\vert \U(\Z/p^l\Z):\pi_l(H) \vert=
\frac{(p-1)p^{\text{min}(L_p(q),L_p(r))-1}}{[\order_p(q),\order_p(r)]}.
$$
We could carry out a similar analysis when $\vert F\vert>1$, though it
would not be so easy to work out some of the indices explicitly.
\end{remark}

\begin{remark}
Theorem~\ref{comp_series} implies in particular that
$C(H)\rtimes_\sigma \Z^2$ satisfies the hypotheses of the
classification theorem of Elliott for AT-algebras
\cite[Theorem~3.2.6]{Ro}. We can tell from the computation of
$K$-theory in Theorem~\ref{two_primes_on_one} that $C(H)\rtimes_\sigma 
\Z^2$ is not a Bunce-Deddens algebra, but it is still closely related
to an odometer. The homeomorphism of $\prod_{k \geq 0}X_k$ onto $H_q$
identifies the action of the first copy of $\Z$ (multiplication by
$q$ on $H_q$) with an odometer (addition of  $1$ with carry over). The
action of the
second copy of $\Z$ (multiplication by $r$ on $H_q$) also acts as a
permutation on each $X_k$: it moves $X_0$ around in a different order,
and this action carries over into $X_1$ when the marker in $X_0$ returns
to the starting point. So we can think of the action of $\Z^2$ as two
independent odometers on the same set. We can normalise the scale so
that either copy of $\Z$ acts by addition of $1$ with carry over, but
not so that both act this way at once.
\end{remark}

\section{The Bost-Connes algebra}\label{BCalgebra}

The Hecke $C^*$-algebra $\mathcal{C}_\Q$ of Bost and Connes \cite{BC}
is
isomorphic to the semigroup crossed product $C^*(\Q/\Z)\rtimes_\alpha
\N^*$. The Fourier transform takes $C^*(\Q/\Z)$
onto the algebra of continuous functions on the compact group
$\mathcal{Z}:=\prod_{p\in \mathcal{P}}\Zp$ and carries $\alpha$ into
the action given by (see \cite[\S3.1]{L}) 
\[
\alpha_n(f)(x)=\begin{cases}f(x/n)&\text{ if }n\text{ divides }x\\ 
0&\text{ otherwise.}\end{cases}
\]

Lemma~\ref{ext_ideals_of_BC_algebra} is valid with $F$ replaced by
$\mathcal{P}$ and $\Z_F$ by $\z$. Thus for $S\subset \mathcal{P}$, 
an application of \cite[Theorem 1.7]{Lar} gives that
$J_S:=C_0(\mathcal{Z}\setminus \z_S)\rtimes_\alpha 
\N^*$ is an ideal of
$\mathcal{C}_\Q=C(\mathcal{Z})\rtimes \N^*$, with quotient isomorphic to
$C(\z_S)\rtimes \N^*$. 
Choose $a\in \mathcal{Z}$ such that $a_p=0 \iff p\in S$. Then
$\{\Q_+^*a\cap \mathcal{Z}\}$ has closure $\mathcal{Z}_S$, so
$C_0(\mathcal{Z}\setminus \mathcal{Z}_S)$ is the kernel of the representation
$\pi_a$ considered in \cite[page~440]{LR3}, and it follows from
\cite[Lemma~4.2]{LPR} that $J_S$ is the kernel of the representation
$\pi_a\times V$ described in \cite[page~440]{LR3}. We can now deduce that $S\mapsto J_S$, as $S$ runs
through the proper subsets of $\mathcal{P}$, is the
parametrisation of $(\Prim C_\Q) \setminus \widehat{\Q_+^*}$ given in
\cite[Theorem~2.8]{LR3}.

\begin{theorem}\label{cofinite_subquotient}
Suppose that $S$ is a proper subset of $\mathcal{P}$.
\begin{enumerate}
\smallskip
\item[(a)] If $\mathcal{P} \setminus S$ is infinite, then
  $J_S=\bigcap_{p\notin S}J_{S\cup\{p\}}.$
\smallskip
\item[(b)] If $0<\vert\mathcal{P} \setminus S\vert<\infty$, then
  \[
\big(\textstyle{\bigcap_{p\notin
    S}}J_{S\cup\{p\}}\big)/J_S\cong\big( C(\U(\Z_{\mathcal{P}\setminus
    S}))\rtimes_\sigma \Z^S \big) \otimes
\K(l^2(\N^{\mathcal{P}\setminus
    S})).
\] 
\item[(c)] $\mathcal{C}_\Q/J_\mathcal{P}$ is isomorphic to
  $C^*(\Q_+^*)=C(\widehat{\Q_+^*})$.
\smallskip
\end{enumerate}
Moreover, $C(\U(\Z_{\mathcal{P}\setminus S}))\rtimes_\sigma \Z^S$ is a
  simple AT-algebra with real rank zero and a unique
  tracial state.
\end{theorem}

It follows from \cite[\S 2]{LR3} that every basic open neighbourhood of
$J_S$ has the form 
\[
U_G=\{J_T\mid T\subset \mathcal{P}, T\cap G=\emptyset\}
\]
for some finite subset $G$ of $\mathcal{P}\setminus S$.
When $\mathcal{P}\setminus S$ is
infinite, there are always lots of $J_{S\cup \{p\}}$ in $U_G$, and thus
$J_S\in \overline{\{J_{S\cup \{p\}}\mid p\notin S\}}$; this says
precisely that $\bigcap_{p\notin S}J_{S\cup\{p\}}\subset J_S$. The
other inclusion is trivial, and (a) follows. Part (c) is true because
$\z_{\mathcal{P}}=\{0\}$. To prove (b) we just need to replace $F$ by
$\mathcal{P}$ and $\Z_F$ by $\z$ in the proof of
Lemma~\ref{identify_small_subquotient_as_scp}.

It remains to prove the statements about
$C(\U(\Z_{\mathcal{P}\setminus S}))\rtimes_\sigma \Z^S$. The Chinese
Remainder Theorem implies that $\Z$ is dense in 
 $\Z_{\mathcal{P}\setminus S}$, and hence $\Z^S=\Z \cap
 \U(\Z_{\mathcal{P}\setminus S})$ is dense in
 $\U(\Z_{\mathcal{P}\setminus S})$. Thus $\Z^S$ acts minimally and
 freely on $\U(\Z_{\mathcal{P}\setminus S})$, and  
$C(\U(\Z_{\mathcal{P}\setminus S}))\rtimes_\sigma \Z^S$ is
simple. However, since $\vert S\vert =\infty$, we cannot apply the
results of \cite{P} to conclude that
$C(\U(\Z_{\mathcal{P}\setminus S}))\rtimes_\sigma \Z^S$ has real rank
zero and stable rank one, as we did in
Section~\ref{finitely_many_primes} for $C(H)\rtimes_\sigma
\Z^S$. Instead we aim to use Theorems 1 and 2 of \cite{BDR}, and to do
this we need to show that $C(\U(\Z_{\mathcal{P}\setminus
  S}))\rtimes_\sigma \Z^S$ is an AH-algebra with the extra property of
slow dimension growth.

Since $\mathcal{P}\setminus S$ is finite, we have as in
(\ref{units_as_inverse_limit}) that 
$\U(\Z_{\mathcal{P}\setminus S})$ is the inverse limit of
the finite groups $\U\big(\Z/\big(\prod_{p\in \mathcal{P}\setminus S}
p^{l_p}\big)\Z\big)$ over $ l=(l_p)\in \N^{\mathcal{P}\setminus S}$. Hence
$\U(\Z_{\mathcal{P}\setminus S})=\varprojlim F_n$, where
\begin{equation}
F_n:= \U\big(\Z/\big(\textstyle{\prod_{p\in
    \mathcal{P}\setminus S}} p^{n}\big)\Z\big). 
 \label{finite_groups_Fl}
\end{equation}
We denote the canonical surjection of
$\U(\Z_{\mathcal{P}\setminus S})$ onto $F_n$ by $\pi_n$. The analogue of
Lemma~\ref{simple_bit_as_union_of_Zd_cp} for $F_n$ and the canonical
action of $\Z^S$ by multiplication on $F_n$ implies that
$C(\U(\Z_{\mathcal{P}\setminus S}))\rtimes_\sigma
\Z^S$ is the closed union of $C^*$-subalgebras isomorphic to 
 $C(F_n)\rtimes \Z^S$.

Towards applying Proposition~\ref{unwinding_Zn_cp}, we note that  
the infinite direct sum $\Z^S$ is the union of the subgroups
$\Z^E$ associated to finite subsets $E$ of $S$. Thus by using an  
argument similar to that in 
Lemma~\ref{simple_bit_as_union_of_Zd_cp} we have that 
$C(F_n)\rtimes \Z^S$ is the closed union of  subalgebras
isomorphic to $C(F_n)\rtimes \Z^E$. However, by choosing a particular
sequence $E_n$ of finite subsets of $S$, we can show that 
$C(\U(\Z_{\mathcal{P}\setminus S}))\rtimes_\sigma \Z^S$ has slow
dimension growth.

Indeed, since $\Z/\big(\prod_p p^n\big)\Z\cong \prod_p
\Z/p^n\Z$, we have $F_n\cong \prod_{p\in \mathcal{P}\setminus S} 
\U (\Z/p^n\Z)$. Thus $F_n$ is a
product of at most $\vert \mathcal{P}\setminus S \vert +1$ cyclic
groups (the $+1$ allows for the possibility that $2\in
\mathcal{P}\setminus S)$,
 and hence has a generating set $\{x_{n, i}
\}$ with at most 
$\vert \mathcal{P}\setminus S \vert +1$ elements. By Dirichlet's
Theorem, there are primes $q_{n, i}$ such that 
\[q_{n, i}\equiv x_{n, i}\pmod{\textstyle{\prod_{p\in
\mathcal{P}\setminus S}}p^n},
\]
and  each $q_{n, i}$ belongs to  $S$ because it is a
unit modulo $\prod_{p\in \mathcal{P}\setminus S}p^n$. Now let
$E_n':=\{q_{n, i}\}$, list the
primes in $S$ as $\{r_n\mid n\in \N\}$, and take 
\[
E_n:=\big(\textstyle{\bigcup_{m\leq n}}E_m'\big)\cup \{r_1, \dots, r_n\}. 
\]
We then have $\pi_n(\Z^{E_n})=F_n$, $E_m\subset E_n$ for $m\leq n$,
and $\bigcup E_n=S$.

We have now realised $C(\U(\Z_{\mathcal{P}\setminus
  S}))\rtimes_\sigma \Z^S$ as the closure of an increasing union
$\bigcup_{n\in \N}B_n$ in which $B_n$ is isomorphic to the
crossed product $C(F_n)\rtimes \Z^{E_n}$ by a transitive action of
$\Z^{E_n}$. By an argument identical to the one at the end of 
Section~\ref{finitely_many_primes} we conclude  that
$C(\U(\Z_{\mathcal{P}\setminus S}))\rtimes \Z^S$ has a unique tracial
state. 

We prove next that 
$C(\U(\Z_{\mathcal{P}\setminus S}))\rtimes \Z^S$ is an AH-algebra
with real rank zero. Proposition
\ref{unwinding_Zn_cp} implies that $B_n\cong C(F_n)\rtimes \Z^{E_n}\cong
C(\T^{\vert E_n \vert}, M_{\vert F_n \vert}(\C))$. But 
\[
\frac {\vert E_n \vert}{\vert F_n \vert}\leq \frac{n(\vert
  \mathcal{P}\setminus S\vert+2)}{\prod_p (p-1)p^{n-1}}\rightarrow 0 \text{
  as }n\rightarrow \infty,
\]
and thus the sequence $B_n\cong C(F_n)\rtimes \Z^{E_n}$ of subalgebras
of $C(\U(\Z_{\mathcal{P}\setminus S}))\rtimes \Z^S$ has slow dimension
growth. It now follows from \cite[Theorem 1]{BDR} that
$C(\U(\Z_{\mathcal{P}\setminus S}))\rtimes \Z^S$ has topological stable
rank one. Since the
projections in $C(\U(\Z_{\mathcal{P}\setminus S}))\rtimes \Z^S$
trivially separate the unique tracial state, \cite[Theorem 2]{BDR} implies
that $C(\U(\Z_{\mathcal{P}\setminus S}))\rtimes \Z^S$ has real rank
zero. The K-groups of $C(\U(\Z_{\mathcal{P}\setminus S}))\rtimes \Z^S$ 
are inductive limits of torsion-free groups, and hence are themselves torsion-free, 
so it follows as in  
Section~\ref{finitely_many_primes} that  
$C(\U(\Z_{\mathcal{P}\setminus S}))\rtimes \Z^S$ is an AT-algebra.

This completes the proof of Theorem~\ref{cofinite_subquotient}.

\appendix
\section{The orders of a prime in groups of units}

For $p$ prime and $m\in \N$ such that $(m, p)=1$, 
we denote by  $\order_{p^l}(m)$ the order of $m$ in
$\U(\Z/p^{l}\Z)$.  
It was shown in \cite[Theorem 3.1]{LPR} that if $p$ is odd, there is a positive
 integer $L_p(m)$ such that
\begin{equation}
\order_{p^l}(m)=\begin{cases}\order_p(m)&\text{ if }1\leq l\leq L_p(m)\\
p^{l-L_p(m)}\order_p(m)&\text{ if }l>L_p(m); \end{cases}
\label{order_of_m_mod_pl}
\end{equation}
the proof uses that the groups $\U(\Z/p^l\Z)$ are cyclic.
We will now show how to
modify the arguments of \cite[\S 3]{LPR} to obtain an analogue of
(\ref{order_of_m_mod_pl}) for $p=2$, in which case $\U(\Z/2^l\Z)$ are no longer
cyclic.

\begin{prop}\label{p2} 
If $m$ is an odd integer and  $m\equiv1\pmod{4}$, then there exists a positive integer 
$K=L_2(m)$ such that
\begin{equation}\label{e7}
\order_{2^l}(m) =
\begin{cases}
    1 &\text{if $1\le l\le K$}\\ 2^{l-K} &\text{if $l>K$;}
\end{cases}
\end{equation}
if $m\equiv3\pmod{4}$, then there exists a positive integer 
$L=L_2(m)$ such that
\begin{equation}\label{e8}
\order_{2^l}(m) =
\begin{cases}
1 &\text{if $l=1$}\\ 2 &\text{if $1<l\le L$}\\ 2^{l-(L-1)} &\text{if $l>L.$}
\end{cases}
\end{equation}
\end{prop}

To prove Proposition \ref{p2} we use general 
properties of cyclic groups as in \cite[\S3]{LPR}. We begin with
a lemma.

\begin{lemma}\label{l4}
Suppose $l\geq 3$. Then the group $\{n\in\U(\Z/2^l\Z) \mid n\equiv1\pmod{4}\}$ is the cyclic 
subgroup ${\langle5\rangle}_l$ of $\U(\Z/2^l\Z)$ generated by $5$.
\end{lemma}

\begin{proof}
Theorem~$2'$ of \cite[Chapter~4.1]{IR} says that
$\vert{\langle5\rangle}_l\vert=2^{l-2}$. For $k\geq 0$ we have
\begin{align*}
5^k ={(4+1)}^k & =
\sum_{n=0}^k \binom kn 4^n = 4\sum_{n=1}^k \binom kn 4^{n-1}+1,
\end{align*}
so $5^k\equiv1\pmod{4}$.
Hence, if $n\equiv 5^k\pmod{2^l}$ for some $0\leq k<2^{l-2}$, then 
$n\equiv1\pmod{4}$. Since the order of $\{n\in\U(\Z/2^l\Z) \mid  
n\equiv1\pmod{4}\}$ is also $2^{l-2}$, the result follows.
\end{proof}

\begin{cor}\label{c1}
An element of $\U(\Z/2^l\Z)$ is congruent to $3\pmod{4}$ 
if and only if it is congruent to ${-5}^k\pmod{2^l}$ for some $k$
satisfying $0\leq k<2^{l-2}$.
\end{cor}

\begin{cor}\label{c2}
Suppose $m\in\Z$ satisfies $m\equiv1\pmod{4}$. Then for every $l>0$ we have
\begin{equation}\label{o2tol}
\order_{2^l}(m)=
\begin{cases}
\order_{2^{l+1}}(m) &\text{if $2$ does not divide $\order_{2^{l+1}}(m)$}\\ 
\order_{2^{l+1}}(m)/2 &\text{if $2$  divides $\order_{2^{l+1}}(m)$.}
\end{cases}
\end{equation}
\end{cor}

\begin{proof}
Since a number is coprime to $2^l$ if and only if it 
is coprime to $2^{l+1}$, the reduction map $\pi :
\U(\Z/2^{l+1}\Z)\rightarrow \U(\Z/2^l\Z)$ is a surjective
homomorphism. Lemma \ref{l4} implies that $m\equiv5^r\pmod{2^{l+1}}$, where 
$r=\order_{2^{l+1}}(5)/\order_{2^{l+1}}(m)=2^{l-1}/\order_{2^{l+1}}(m)$. Hence, by  
applying \cite[Lemma 3.2]{LPR} to the restriction of 
$\pi$ to a homomorphism of 
${\langle 5\rangle}_{l+1}$ onto ${\langle 5\rangle}_l$, we have
\begin{align*}
\order_{2^l}(m)&=\order(\pi(5^r))\\
 & =\begin{cases}
2^{l-1}/(2^{l-1}/\order_{2^{l+1}}(m),2^{l-1}) &\text{if $2^{l-1}$ divides 
$2^{l-1}/\order_{2^{l+1}}(m)$}\\ 
2^{l-1}/\big(2(2^{l-1}/\order_{2^{l+1}}(m),2^{l-1})\big) &\text{if $2^{l-1}$ does 
not divide $2^{l-1}/\order_{2^{l+1}}(m)$,}
\end{cases}
\end{align*}
which simplifies to (\ref{o2tol}).
\end{proof}

\begin{proof}[Proof of Proposition \ref{p2}]
Suppose first that  $m\equiv1\pmod{4}$. For fixed $N$, there 
exists an $l\in\N$ satisfying $m^N<2^l$. Then $o_{2^l}(m)>N$ and 
hence the sequence $\{\order_{2^l}(m) \mid l\in\N\}$ must be unbounded. In 
particular, $\{\order_{2^l}(m)\}$ is not a constant sequence. Let $K$ be 
the first integer such that $\order_{2^K}(m)<\order_{2^{K+1}}(m)$. Then 
$\order_{2^l}(m)=\order_2(m)=1$ for $1\le l\le K$, and by Corollary
\ref{c2} we have  
$\order_{2^{K+1}}(m)=2\order_2(m)=2$. Since $\order_{2^{K+1}}(m)$ divides 
$\order_{2^l}(m)$ for all $l>K$, it follows that $2$ divides
$\order_{2^l}(m)$  
for all $l>K$. We now apply Corollary \ref{c2} $l-K$ times to 
deduce that $\order_{2^l}(m)=2^{l-K}\order_{2^K}(m)=2^{l-K}$.

Now suppose that $m\equiv 3\pmod{4}$. Certainly  $\order_2(m)=1$. For
$l>1$, Corollary \ref{c1} tells us that $m\equiv -5^k$ for some $0\leq
k<2^{l-2}$. Thus $m^2\equiv  5^{2k}\pmod{2^l}$, and therefore $m^2\in
{\langle 5\rangle}_l$. Let $L$ be the first integer such that
$\order_{2^L}(m^2)<\order_{2^{L+1}}(m^2)$. Applying Corollary~\ref{c2} to 
$m^2$ and repeating the argument of the preceding paragraph gives
(\ref{e8})  because $\order_{2^l}(m)=2\order_{2^l}(m^2)$.
\end{proof}

We now need to extend these results to cover actions on $\U(\Z_F)$ for
an arbitrary finite set $F$ of primes. We write $F=\{p_1, \dots
,p_{n}\}$ and fix a prime $q$ which is not in $F$. We denote by 
$\order_{(l_1,\dots,l_n)}(q)$ the order of $(q,\dots,q)$ in 
$\prod_{i=1}^n\U(\Z/{p_i}^{l_i}\Z)$.

\begin{prop}\label{nothyfacts}
There exist positive integers  
$K_1,\dots,K_n$ and $d$ such that
\begin{equation}\label{order_mod_nprimes}
\order_{(K_1+l_1,\dots,K_n+l_n)}(q)=d{p_1}^{l_1}\dots{p_n}^{l_n}
\end{equation}
for every $(l_1,\dots ,l_n)\in \N^F.$
\end{prop}

\begin{proof}
Suppose first that $p_1,\dots,p_n$
 are distinct odd primes, and let $L_{p_i}(q)$ be as in (\ref{order_of_m_mod_pl}). Let 
\[
z_i:=\max\{z\mid {p_i}^z\text{ divides } \order_{p_j}(q) \text{ for
  some }j\in\{1,\dots,n\}\},
\] 
and define 
$K_i:= L_{p_i}(q) + z_i$ and $d:=[\order_{p_1}(q),\dots,\order_{p_n}(q)]$, 
where $[r_1, \dots , r_n]$ is the least common multiple of the integers $r_i$. 
In general,  if $g_i$ are elements of order $r_i$ in finite groups  
$G_i$, then the order of $(g_1, \dots ,g_n)$ in
$G_1\times \dots \times G_n$ is $[r_1, \dots , r_n]$.
Thus from the properties of $L_{p_i}(q)$  we obtain
\begin{align}
\label{order_of_q_for_p_i_odd}
\order_{(K_1+l_1,\dots,K_n+l_n)}(q) & =
[p_1^{(K_1+l_1)-L_{p_1}(q)}\order_{p_1}(q), \dots ,
p_n^{(K_n+l_n)-L_{p_n}(q)}\order_{p_n}(q)] \\
&=[p_1^{z_1+l_1}\order_{p_1}(q), \dots ,
p_n^{z_n+l_n}\order_{p_n}(q)] \notag \\
&={p_1}^{l_1}\dots{p_n}^{l_n}[\order_{p_1}(q),\dots,\order_{p_n}(q)]\notag,
\end{align}
which is (\ref{order_mod_nprimes}).

Now suppose that $2\in F$, say $p_1=2$. If $q\equiv1\pmod{4}$, we let 
\[
z_i:=\max\{z\mid {p_i}^z\text{ divides }
\order_{p_j}(q)\text{ for some } j\in\{2,\dots,n\}\},
\]
and define $K_1:=L_2(q)+z_1$, $K_i:=L_{p_i}(q)+z_i$ for
$i>1$, and $d:=[\order_{p_2}(q),\dots,\order_{p_n}(q)]$. 
Reasoning as in
(\ref{order_of_q_for_p_i_odd}) gives (\ref{order_mod_nprimes}). 

If $q\equiv3\pmod{4}$, we let 
\begin{align*}
z_1&=\max(1,\max\{z\mid 2^z\text{ divides }
\order_{p_j}(q) \text{ for some }j\in\{2,\dots,n\}\}),\\
z_i&=\max\{z\mid {p_i}^z\text{ divides }
\order_{p_j}(q) \text{ for some }j\in\{2,\dots,n\}\} 
\end{align*} 
for $i>1$, and define $K_1:=L_2(q)+z_1-1$, $K_i:=L_{p_i}(q)+z_i$ for
$i>1$ and $d:=[2, \order_{p_2}(q),\dots,\order_{p_n}(q)]$. 
Again, reasoning as in (\ref{order_of_q_for_p_i_odd}) 
gives  (\ref{order_mod_nprimes}). 
\end{proof}

\begin{cor}\label{indexform}
The closure $H$ of $q^\Z$ in $\U(\Z_F)$ is a subgroup of finite index 
\[
\vert \U(\Z_F):H\vert=\big(\textstyle{\prod_{i=1}^n (p_i-1)p_i^{K_i-1}}\big)/d.
\]
\end{cor}

\begin{proof}
Apply Proposition~\ref{nothyfacts} to ${\bf l}=(l,l,\dots,l)$ to see
that
$\vert \pi_l(H)\vert=d\big(\textstyle{\prod_{i=1}^n p_i^{l}}\big)$
for large $l$, and the result follows from \cite[Lemma~3.7]{LPR}.
\end{proof}

\end{document}